\documentclass[a4]{article}

\usepackage[utf8]{inputenc}
\usepackage{mathtools}
\usepackage{placeins}
\usepackage{subfig}
\usepackage{amsmath}
\usepackage{amsfonts}
\usepackage{amssymb}
\usepackage{url}
\usepackage{tikz}
\usetikzlibrary{shapes,arrows}
\usetikzlibrary{arrows.meta}
\usetikzlibrary{calc}
\usetikzlibrary{positioning}
\usepackage{pgfplots}
\usepackage{pgfplotstable}
\pgfplotscreateplotcyclelist{mycycle}{%
	solid, every mark/.append style={solid, fill=gray}, mark=*\\%
	dashed, every mark/.append style={solid, fill=gray}, mark=square*\\%
	dotted, every mark/.append style={solid, fill=gray}, line width=0.8pt, mark=diamond*\\%
	dashdotted, every mark/.append style={solid, fill=gray, rotate=180},mark=triangle*\\%	
	densely dashed, every mark/.append style={solid, fill=gray}, mark=otimes*\\%
	loosely dashed, every mark/.append style={solid, fill=gray}, mark=triangle*\\%
	loosely dotted, every mark/.append style={solid, fill=gray}, mark=oplus*\\%
	densely dotted, every mark/.append style={solid, fill=gray}, mark=pentagon*\\%
}
\usepgfplotslibrary{groupplots}
\pgfplotsset{
	compat=newest,
	grid=major,
	cycle list name=mycycle,
	yminorticks=false,
	tick align=outside,
	tick pos=left,
	legend pos=outer north east,
	legend cell align=left,
	height=0.2\textwidth,
	width=0.2\textwidth,
	every axis/.append style={scale only axis},
	legend style={font=\small},
	label style={font=\small},
	tick label style={font=\small},
}

\usepgfplotslibrary{colormaps}
\pgfplotsset{%
	/pgfplots/colormap={visit}{%
		rgb255(0cm)=(0,0,255)
		rgb255(1cm)=(0,255,255)
		rgb255(2cm)=(0,255,0)
		rgb255(3cm)=(255,255,0)
		rgb255(4cm)=(255,0,0)
	}
}

\usepackage{amsthm}
\newtheorem{definition}{Definition}
\newcommand{\normVector}[1]{\left| {#1} \right|}
\newcommand{\cutcell}{K_{j, \mathfrak{s}}}

\newcommand{\jump}[1]{\left[\!\left[{{#1}}\right]\!\right]}

\newcommand{\shockpos}{x_\mathrm{s}}
\newcommand{\interfacepos}{x_\mathfrak{I}}
\newcommand{\pzeroindicator}{\mathcal{I}^{P0}}
\newcommand{\rhoindicator}{\mathcal{I}^{\rho}}
\newcommand{\momindicator}{\mathcal{I}^{m}}
\newcommand{\diff}[1]{\,\mathrm{d} #1}

\newcommand{\shockLevelSet}{\varphi_\mathrm{s}}
\newcommand{\shockInterface}{\mathfrak{I}_\mathrm{s}}

\title{An Extended Discontinuous Galerkin Method for High-Order Shock-Fitting}
\author{Markus Geisenhofer, Florian Kummer, and Martin Oberlack}
\begin{document}

\maketitle

% ###############################################################################
% ###############################################################################
\begin{abstract}
We present a sub-cell accurate shock-fitting technique 
using a high-order extended discontinuous Galerkin (XDG) method,
where a computational cell of the background grid is cut into two cut-cells at the shock position. 
Our technique makes use of a sharp interface description where the shock front is implicitly defined by means of the zero iso-contour of a level-set function.
A novel implicit pseudo-time-stepping procedure is employed to correct the position of the shock front inside the cut background cell by using cell-local indicators, since the position and shape of shock waves are not known a priori for the general, multi-dimensional case. 
This iterative correction terminates if the shock front has converged to the exact position. 
The procedure is demonstrated for the test case of a one-dimensional stationary normal shock wave.
Furthermore, the underlying sharp interface approach drastically reduces the complexity of the grid handling, since a simple Cartesian background grid can be employed.
\end{abstract}
% ###############################################################################
% ###############################################################################

% ###############################################################################
% ###############################################################################
\section{Introduction}
% ###############################################################################
% ###############################################################################
Within the past decades, high-order discontinuous Galerkin (DG) methods~\cite{reed1973,cockburn2000,warburton1999} have become a popular tool to discretize systems of partial differential equations (PDE) in the context of fluid flows. DG methods combine advantages of classical finite volume methods (FVM) and finite element methods (FEM), such as conservativity, their local character due to the weak coupling of neighboring cells via cell boundaries, and the applicability to structured and unstructured grids. Furthermore, these properties render DG methods highly attractive for high performance computing applications~\cite{altmann2013a,dryja2016}.

High-order extended discontinuous Galerkin (XDG) methods are a sub-class of DG methods in which the geometry or domain of interest do not conform with the computational domain.
The approximation space is enriched with an additional set of degrees of freedom (DOF) in cells which contain a non-smooth or discontinuous solution of the PDE, for example, as being present in the case of a water bubble in air in the context of multi-phase flows. In extended methods, an interface usually separates the computational domain into two disjoint regions.
The idea of extended methods goes back to the XFEM approach by Mo\"{e}s~\cite{moes1999} who investigated crack propagation in solid mechanics. Adaptions to incompressible multi-phase flows were published, e.g., by Gro{\ss} and Reusken~\cite{gross2007} and Fries~\cite{fries2009}. Bastian and Engwer~\cite{bastian2009} were the first to bring the extended idea to DG methods, focusing on two- and three-dimensional, elliptic, scalar model problems on complex-shaped domains. This and further related works share the idea of using piece-wise planar, triangular sub-cells for the integration of the weak forms in cut-cells. Heimann~\cite{heimann2013} as well as Kummer~\cite{kummer2016} published extensions to incompressible flows based on a (symmetric) interior penalty discretization~\cite{arnold1982}. Kummer additionally employed a cell-agglomeration technique in order to enhance the conditioning of the system matrices, also applying the Hierarchical-Moment-Fitting (HMF) quadrature scheme as proposed by Müller~\cite{muller2013a}. 

Immersed boundary methods (IBM), originally proposed by Peskin~\cite{peskin1972} in the context of blood flows, summarized in detail by Mittal and Iaccarino~\cite{mittal2005}, and later adapted to DG methods by Fidkowski and Darmofal~\cite{fidkowski2007}, can be seen as a sub-class of extended methods in which the geometry representation and the discretization are separated to a large extent. A DG IBM can be considered as an XDG method where one region of interest describes the fluid phase and the other one represents the geometry, which is void. Further works were published by Qin and Krivodonova~\cite{qin2013} who presented a DG IBM using explicit time-integration schemes for smooth solutions of the Euler equations on a Cartesian grid. Müller et al.~\cite{muller2017} extended the work by Qin and Krivodonova~\cite{qin2013} to the compressible Navier-Stokes equations by adapting the HMF quadrature scheme~\cite{muller2013a} for the application in compressible flows in combination with a cell-agglomeration strategy. Recently, Geisenhofer et al.~\cite{geisenhofer2019} combined the DG IBM with a shock-capturing strategy based on artificial viscosity for the application in supersonic compressible flows. Additionally, a local time stepping scheme has been applied for computational speed, since non-agglomerated cut-cells restrict the maximum admissible time-step size drastically.

Discontinuous flow phenomena such as shock waves can arise in supersonic flows. In general, DG methods are susceptible to stability issues which are caused by oscillating polynomial solutions in the vicinity of discontinuities or under-resolved regions of the computational domain~\cite{persson2006,barter2010,klockner2011}. Remedies are classical limiting approaches such as ENO/WENO schemes~\cite{shu1988,shu1989} or, e.g., a posteriori limiting approaches~\cite{giani2014,dumbser2014}. Another class are artificial viscosity methods, which date back to the early work by Von Neumann and Richmyer~\cite{vonneumann1950}. An additional `artificial', second-order, diffusive term smooths discontinuities over a layer~$O (h/P)$, which can be adequately resolved by the numerical scheme. Persson and Peraire~\cite{persson2006} proposed a sensor based on the modal-decay of the local solution to detect troubled cells and a general scaling for determining the amount of artificial viscosity. Further improvements and applications were published by Barter and Darmofal~\cite{barter2010}, Klöckner et al.~\cite{klockner2011}, and Lv et al.~\cite{lv2016}. However, a clear disadvantage of shock-capturing strategies is that low-order perturbations such as acoustic waves are usually damped out entirely. This still remains an issue to tackle in the context of high-order methods~\cite{wang2013}.

By contrast, classical shock-fitting approaches, first published by Emmons \cite{emmons1944,emmons1948} in the 1940's, have in common that one boundary of the computational domain is fitted to the shock front so that only the post-shock region has to be computed. The conditions at the shock are prescribed via the Rankine-Hugoniot conditions. Moretti and Abbett~\cite{moretti1966} were the first to solve the supersonic blunt-body problem with a practical technique in the 1960's. It relies on a mapping of the physical, time-dependent domain onto a rectangular reference domain in the context of a finite difference method (FDM). Additionally, they applied a time-marching technique to compute the steady-state solution. We recommend the work by Salas~\cite{salas2010} for a historical and detailed introduction. Variants of the original approach which are based on Fourier and Chebyshev approximations were developed, e.g., by Zang et al.~\cite{zang1983} and Hussaini et al.~\cite{hussaini1985a}. Combinations were made by Wu and Zhu~\cite{wu1996} and Kopriva~\cite{kopriva1991} in the context of multi-domain methods, where the domain of interest is split into several sub-domains according to an a-priori knowledge about the flow. Then, the discontinuities are fitted to each sub-domain separately. Lately, Romick and Aslam~\cite{romick2017} applied shock-fitting to detonation problems~\cite{henrick2008,short2008}. A computational boundary was fitted to shock fronts and material interfaces, while solving the two-dimensional reactive Euler equations.
Recently, another DG shock-fitting variant was developed by Corrigan et al.~\cite{corrigan2019}, who proposed a moving DG method with interface condition enforcement (MDG-ICE). Their approach relies on a weak formulation which enforces the interface condition separately from the conservation law. The discrete grid geometry is added as an additional variable to the solver. Their approach detects a priori unknown discontinuities via the interface condition enforcement and satisfies the conversation law by moving the computational grid in every time-step. They showed results with high-order accuracy for steady and unsteady test problems such as for a shock tube and a Mach-3 bow shock.
The approach of Zahr and Persson~\cite{zahr2018} is similar but conceptually different in a sense that they state an optimization problem based on the variation of the discrete solution from its element-wise average. Additionally, they incorporate a grid distortion metric to penalize oscillating solutions and distorted grids. Zahr and Persson showed high-order convergence rates for the inviscid Burgers' equation as well as for a transonic flow through a nozzle and a supersonic flow around a cylinder.

To the best of our knowledge, no efforts have been yet undertaken for a full high-order XDG method in the context of shock-fitting. We aim for solving compressible flows with high Mach numbers, where a sharp interface description of the shock front is employed to separate the pre- and post-shock state. In contrast to the previously mentioned works in the context of DG methods, a simple Cartesian background grid is employed, reducing the complexity of the grid handling to a minimum, since the shock front is inherently incorporated into the underlying discretization.

This work is structured as follows: We state a non-dimensional, conservative form of the Euler equation for inviscid compressible flow in Section~\ref{sec:governing_equations}. In Section~\ref{sec:discretization}, we briefly introduce the employed XDG method in a one-dimensional setting.
In Section~\ref{sec:sub-cell_accurate_correction}, we present an implicit pseudo-time-stepping procedure which corrects the shock interface position inside a cut background cell. To achieve this, we apply suitable indicators to determine the shifting direction of the shock interface. The new interface position is obtained by a simple bisection algorithm. We close this work with a summary and an outlook for future work in Section~\ref{sec:conclusion}.

% ###############################################################################
% ###############################################################################
\section{Governing Equations}	
\label{sec:governing_equations}
% ###############################################################################
% ###############################################################################
This work deals with supersonic inviscid compressible flow. 
This kind of flow can be mathematically described by the Euler equations which we state in a non-dimensional, 
conservative form in one spatial dimension as
\begin{align}	\label{eq:eulerEquations}
	\frac{\partial \vec{U}}{\partial t} 
  + \frac{\partial \vec{F}(\vec U)}{\partial x} 
  = 0\,,
\end{align}
where $\vec{U}$ denotes the vector of conserved quantities 
\begin{align}	\label{eq:vectorConservedQuantities}
	\vec{U} =	\left( \begin{array}{c}
		\rho\\
		\rho u  \\
		\rho E
	\end{array} \right)\,,
\end{align}
and $\vec{F}$ is the convective flux given by
\begin{align}	\label{eq:convectiveFluxes}
	\vec{F}_1 =	\frac{1}{\gamma \mathrm{M}_\infty^2}\left( \begin{array}{c}
		\rho u \\
		\rho u^2 + p\\
		u (\rho E + p)
	\end{array} \right) .
\end{align}
In Equations~\eqref{eq:eulerEquations} to~\eqref{eq:convectiveFluxes}, $t$ is the time, 
$x$ is the spatial coordinate, 
$\rho$ is the fluid density, 
$u$ is the velocity, 
$\rho E$ is the total energy, 
$p$ is the pressure, 
$\gamma$ is the heat capacity ratio which is $\gamma = 1.4$ for standard air, 
and $\mathrm{M}_\infty$ is the dimensionless reference Mach number which we set to $\mathrm{M}_\infty = 1 / \sqrt{\gamma}$. 
Consequently, the non-dimensional Euler equations match their equivalent with dimensions~\cite{feistauer2003}, 
which enables the direct use of initial conditions of any test cases without further modifications. 
The total energy consists of the sum of the internal energy~$\rho e$ and the kinetic energy, i.e., $\rho E = \rho e + 1/2\, \rho u^2 $.

The Euler equations have to be closed by a suitable equation of state for the pressure~$p(\rho, e)$. 
In this work, we consider a calorically perfect gas which can be modeled by the ideal gas law
\begin{align}	\label{eq:ideal_gas_law}
	p(\rho, e) = (\gamma - 1) \rho e \,.
\end{align}
The local Mach number is given by the relation~$\mathrm{M} = |u| / a$, 
where $a = \sqrt{\gamma p / \rho}$ denotes the local speed of sound.

% ###############################################################################
% ###############################################################################
\section{Discretization}	
\label{sec:discretization}
% ###############################################################################
% ###############################################################################
In this section, we briefly introduce the spatial discretization of the employed XDG method in one spatial dimension and apply it to the Euler equations~\eqref{eq:eulerEquations} to~\eqref{eq:convectiveFluxes}. In particular, our approach makes use of a sharp interface description by defining interfaces by means of the zero iso-contour of a level-set function. Furthermore, different level-set functions can be applied for the representation of the geometry, e.g., the surface of solid bodies, and the shock front. The interested reader is referred to the works by Kummer et al.~\cite{kummer2020} and Smuda and Kummer~\cite{smuda2020a} for details about the employed XDG method and its features. The presented methodology is implemented in the open-source software package~\emph{BoSSS}\footnote{\url{https://github.com/FDYdarmstadt/BoSSS}, accessed on 11/24/2020}, which features a variety of applications in the context of boundary-fitted and unfitted DG discretizations for incompressible multi-phase flows~\cite{kummer2016,smuda2020a}, compressible flows~\cite{muller2017,geisenhofer2019}, and particulate flows~\cite{krause2017}.

In order to incorporate the sharp interface description into the discretization approach, 
we introduce a simply connected computational domain
$\Omega_h = (x_{\text{left}}, x_{\text{right}}) \subset \mathbb{R}$, 
which we partition into two disjoint regions, or sub-domains $\mathfrak{A}$ and~$\mathfrak{B}$. 
These are separated by the interface~$\mathfrak{I} = \{ x_{\mathfrak{I}} \} $,
which corresponds to a single point in the one-dimensional setting.
Thus, we obtain the partitioning
\begin{equation}
\Omega_h = 
\mathfrak{A} \, \dot \cup \,  \mathfrak{I} \, \dot \cup \, \mathfrak{B}
= (x_{\text{left}}, x_{\mathfrak{I}} )  
   \, \dot \cup \, \{ x_{\mathfrak{I}} \}  
       \, \dot \cup \,  (  x_{\mathfrak{I}} , x_{\text{right}}  )\,.
\label{eq:partioning}
\end{equation}

In the following, we briefly introduce several definitions:
\begin{definition}[Basic notations]
\begin{itemize}
\item
We discretize the computational domain $\Omega_h$ into a discrete set of non-overlapping cells~$\mathfrak{K}_h = \{ K_1, \dots, K_J \}$ 
with $K_j = (x_j, x_{j+1})$ and 
~$x_{\text{left}} = x_1 \leq x_2 \leq \ldots \leq x_{J+1} = x_{\text{right}}$.
The characteristic length scale is given by~$h = \max_{1 \leq j \leq J} ( x_{j+1} - x_{j} )$.

\item
The set of all cell boundaries or cell edges, respecitvely, is given by~$
\Gamma \coloneqq \cup_j \partial K_j \cup \mathfrak{I}
 = \{  x_1, \ldots x_{J}, x_{\mathfrak{I}} \}
$. It can be split into $\Gamma = \Gamma_\mathrm{int} \cup \Gamma_\mathrm{D} \cup \Gamma_\mathrm{N}$, where $\Gamma_\mathrm{int} = \Gamma \setminus \partial \Omega$ denotes the set of all inner edges, and $\Gamma_\mathrm{D}$ and $\Gamma_\mathrm{N}$ denote the sets of boundary edges with Dirichlet and Neumann boundary conditions, respectively.
		
\item
We introduce a field of normal vectors~$n_\Gamma$ on~$\Gamma$,
which is simply set to $n_\Gamma = -1$ at $x_{\text{left}}$ and to $n_\Gamma = 1$ everywhere else.
This is mainly introduced as a general notation for the two- and three-dimensional setting.
\end{itemize}
\end{definition}

\begin{definition}[Cut-cells and cut-cell grid] \label{def:cut-cell_grid} 
On the background grid $\mathfrak{K}_h = \{ K_1,\allowbreak \dots, K_J \}$, we define a cut-cell~$\cutcell$ of some sub-domain~$\mathfrak{s} \in \{ \mathfrak{A}, \mathfrak{B} \}$ as
\begin{align}	\label{eq:cut_cells}
	\cutcell \coloneqq K_j \cap \mathfrak{s}\,.
\end{align}
The set of all cut-cells builds the cut-cell grid
\begin{align}
	\mathfrak{K}_h^\mathrm{X}
	=
	\{ K_{1, \mathfrak{A}}, K_{1, \mathfrak{B}}, \dots,  K_{J, \mathfrak{A}}, K_{J, \mathfrak{B}}\}\,.
\end{align}
A background cell~$K_j$ which is cut by the interface~$\mathfrak{I}$, i.e., $\oint_{K_j \cap \mathfrak{I}} 1 \diff{S} \neq 0$, 
is divided into two cut-cells~$K_{J, \mathfrak{A}}$ and $K_{J, \mathfrak{B}}$. 
The cell~$K_j$ from the background grid is recovered if it is entirely contained in one sub-domain, 
i.e., $K_{J, \mathfrak{A}} = K_j$ with $K_{J, \mathfrak{B}} = \emptyset$, or $K_{J, \mathfrak{B}} = K_j$ with $K_{J, \mathfrak{A}} = \emptyset$.
\end{definition}
In general, the presented XDG method can be seen as a DG method applied on a cut-cell grid. This motivates the definition of the corresponding polynomial space in the subsequent definition.
\begin{definition}[Extended discontinuous Galerkin (XDG) space]
	We define the broken polynomial, cut-cell XDG space
	\begin{align}	\label{eq:xdg_space}
		\begin{split}
		\mathbb{P}_P^\mathrm{X} (\mathfrak{K}_h)
		\coloneqq
		\{
		f \in L^2 (\Omega); \,
		\forall K \in \mathfrak{K}_h:
		f\vert_{K \cap \mathfrak{s}} \textrm{ is polynomial and }
		\mathrm{ deg} (f\vert_{K \cap \mathfrak{s}} \leq P)
		\}\\
		=
		\mathbb{P}_P (\mathfrak{K}_h^\mathrm{X})
		\end{split}
	\end{align}
	with a total polynomial degree~$P$. 
\end{definition}
\begin{definition}[Jump operator]
	Last, we define a jump operator~$\jump{\psi}$ which acts on the set of edges~$\Gamma$ as
	\begin{align}
		\jump{\psi} \coloneqq
		\begin{cases}
			\psi^- - \psi^+\,,	&\qquad x \in \Gamma_\mathrm{int}\, ,\\
			\psi^-\,, 			&\qquad x \in \partial \Omega\, ,\\
		\end{cases}
	\end{align}
    where $\psi^-$ and $\psi^+$ denote the limit of $\psi$ at $x$ when approaching from $-n_{\Gamma}$ and $+n_{\Gamma}$, respectively.
\end{definition}
Obviously, for an arbitrarily placed interface, the size of some cut-cell may be arbitrarily small. 
Especially when one considers higher spatial dimensions, 
small and ill-shaped cut-cells can significantly limit the global time-step size 
in the context of explicit time-stepping schemes and increase the condition numbers of the cell-local mass matrices~\cite{muller2017}. 
A potential remedy is the application of cell-agglomeration techniques, 
which remove undesired cut-cells from the computational grid. 
We refer to the works by Kummer et al.~\cite{kummer2020} and Müller et al.~\cite{muller2017} for further details and implementation issues.
For this work, however, \emph{we assume that the interface is well-placed}:
This means that the volume fraction of a cut-cell is above a certain threshold, 
i.e. $ \normVector{\cutcell} / \normVector{K_j} \geq  \delta_\mathrm{agg}$ with 
the threshold set to $\delta_\mathrm{agg} = 0.3$.

The convective fluxes of the Euler Equations~\eqref{eq:eulerEquations} to \eqref{eq:convectiveFluxes} are discretized by means of an approximate Riemann solver based on Godunov's method following Section 4.9 of the textbook~\cite{toro2009} by Toro. Note that classical choices such as the HLLC flux~\cite{toro1994} may fail due to inappropriate wave speed estimates in cases where high-sped flows impinge on solid stationary walls~\cite{toro2009}. In addition to such configurations, we could confirm this effect for pseudo-two-dimensional computations of a stationary normal shock wave, which is located on a cell boundary or an interface. In such cases, round-off errors determine the sign of several terms during the flux evaluation and, thus, prohibit a robust and reliable application of the HLLC flux.

% ###############################################################################
% ###############################################################################
\section{Sub-Cell Accurate Correction of the Shock Interface Position}	
\label{sec:sub-cell_accurate_correction}
% ###############################################################################
% ###############################################################################
In this section, we present a novel algorithm for the correction of the shock interface position~$\interfacepos$ in a cut background cell, i.e., $\interfacepos \in K_j$ with~$\oint_{K_j \cap \mathfrak{I}} 1 \diff{S} \neq 0$, see also Definition~\ref{def:cut-cell_grid}. To this end, we assume that a cell-accurate guess for~$\interfacepos$ is already known.

In Section~\ref{sec:indicators}, we introduce three indicators, which are based on a zeroth-order projection of the local solution and the normal shock relations, respectively, in order to correct the shock interface position~$\interfacepos$. In Section~\ref{sec:implicit_pseudo_ts}, a novel implicit pseudo-time-stepping procedure is derived. This correction procedure is applied until the shock interface position~$\interfacepos$ has converged to the exact shock position~$\shockpos$, i.e., $\lim_{l \rightarrow \infty} x_{\mathfrak{I}}^l = \shockpos$. After the termination of this procedure, the entire flow field can be advanced in time like in a standard DG computation. In particular, we show the results of a one-dimensional proof of concept for a stationary normal shock wave in Section~\ref{sec:implicit_pseudo_ts_proof}.

\subsection{Indicators}	\label{sec:indicators}
Subsequently, we present three different indicators in order to determine the shifting direction of the shock interface~$\mathfrak{I}$ towards the exact shock position~$\shockpos$ inside a cut background cell~$K_j$.

\subsubsection{\textit{P0}-Indicator}	\label{sec:p0_indicator}
First, we present an indicator which is based on a zeroth-order projection of the local solution in both cut-cells~$\cutcell \in \mathfrak{K}_h^\mathrm{X}$ with~$\cutcell \neq \emptyset$ and~$\mathfrak{s} \in \{\mathfrak{A}, \mathfrak{B}\}$. We denote the local solution by
\begin{align}	\label{eq:cell_local_sol}
	\psi_{j, \mathfrak{s}} (x, t)
	\coloneqq
	\psi\vert_{\cutcell} \in \mathbb{P}_P^\mathrm{X} (\{ \cutcell \})\,.
\end{align}
If the shock interface position~$\interfacepos$ does not coincide with the exact shock position~$\shockpos$, i.e., $\interfacepos \neq \shockpos$, the polynomial solution oscillates so that it is hardly possible to extract reasonable physical information from it. Thus, we remove the high-order modes by solely extracting the zeroth-order mode
\begin{align}	\label{eq:zeroth_order_mode}
	\psi_{j, \mathfrak{s}}^{P0}
	=
	\Pi^{P0} (\psi_{j, \mathfrak{s}})\,,
\end{align}
using the projection operator
\begin{equation} \label{eq:p0_proj_op}
	\begin{aligned}	
		\Pi^{P0} :
		\mathbb{P}_P^\mathrm{X} (\{ \cutcell \})
		&\rightarrow
		\mathbb{P}_{P=0}^\mathrm{X} (\{ \cutcell \})\,,
		\\
		\psi &\mapsto \psi^{P0}\,,
	\end{aligned}
\end{equation}
with the essential property~$\langle \psi - \psi^{P0}, \vartheta \rangle = 0$, $\forall \vartheta \in \mathbb{P}_{P=0}^\mathrm{X} (\{ \cutcell \})$. 
The $L^2$ scalar product is denoted by~$\langle \cdot \,, \cdot \rangle$. 
Exemplarily, we show the numerical solution for a stationary normal shock wave with 
a Mach number of~$\mathrm{M}_\mathrm{s} = 1.5$ 
located at~$\shockpos = 0.55$ in a cut background cell~$K_j = [0.5, 0.6]$ 
with the sub-domains~$\mathfrak{A}$ and~$\mathfrak{B}$ in Figure~\ref{fig:p0_indicator}.
\begin{figure}[tbp]
	\centering
	\begin{tikzpicture}
		\begin{axis}[
			xmin=0.45, xmax=0.65,
			ymin=0.8, ymax=2.2,
			xlabel=$x$,
			ylabel=Density~$\rho$,
			xmajorgrids=true,
			ymajorgrids=false,
			no markers,
			xtick={0.5, 0.55, 0.6},
			xticklabels={$0.5$, $\shockpos$, $0.6$},
			width=0.27\textwidth,
			height=0.27\textwidth
			]
			
			\draw [fill=lightgray,draw=none] (0.45,0.8) rectangle (0.5,2.2);
			\draw [fill=lightgray,draw=none] (0.6,0.8) rectangle (0.65,2.2);
			
			\pgfplotstableread{XDG_SSW_Initial_Condition.txt}\datafile
			\addplot+[color=black, dashed] table[skip first n=1] {\datafile};
			\addlegendentry{Shock}
			
			\pgfplotstableread{1_rho.txt}\datafile
			\addplot+[solid] table[skip first n=1] {\datafile};
			\addlegendentry{$P=2$}
			
			\pgfplotstableread{1_rho_p0.txt}\datafile
			\addplot+[dashdotted] table[skip first n=1] {\datafile};
			\addlegendentry{$P0$-projection}
			
			\draw[color=red] (axis cs:0.57,0.8) -- (axis cs:0.57,2.2);
			\node[color=red, anchor=south west] at (axis cs:0.57,0.8) {\small$\interfacepos$};	
			
			\node[color=black, anchor=south west] at (axis cs:0.451,0.77) {\small$\rho_\mathrm{pre}$};	
			\node[color=black, anchor=south west] at (axis cs:0.6,1.65) {\small$\rho_\mathrm{post}$};	
			
			\node[color=black, anchor=south west] at (axis cs:0.505,1.97) {\small$\mathfrak{A}$};	
			\node[color=black, anchor=south west] at (axis cs:0.594,1.97) {\small$\mathfrak{B}$};	
			
		\end{axis}	
	\end{tikzpicture}
\caption{$P0$-projection of the local solution in a cut background cell~$K_j = [0.5,0.6]$.
A stationary shock wave with a Mach number of $\mathrm{M}_\mathrm{s} = 1.5$ is located at $\shockpos = 0.55$. The solution oscillates, since the shock interface position~$\interfacepos$ does not coincide with the exact shock position~$\shockpos$, i.e., $\interfacepos \neq \shockpos$.
The $P0$-projection considers only the zeroth-order modes in sub-domain~$\mathfrak{A}$ and $\mathfrak{B}$, respectively. The exact density values~$\rho_\mathrm{pre}$ and $\rho_\mathrm{post}$ can be extracted from the near-band~$\mathfrak{K}_h^\mathrm{near}$ (light gray) for later usage.
}
\label{fig:p0_indicator}
\end{figure}
The high-order solution ($P=2$) oscillates, since the shock interface position and the exact shock position do not coincide, i.e., $\interfacepos \neq \shockpos$.
Additionally, we show the corresponding $P0$-projections in both sub-domains~$\mathfrak{A}$ and~$\mathfrak{B}$. 
In general, the projections in~$\mathfrak{A}$ and~$\mathfrak{B}$ are different, 
i.e., $\psi_{j, \mathfrak{A}}^{P0} \neq \psi_{j, \mathfrak{B}}^{P0}$, since we have a separated set of DOF for each sub-domain.

We construct an indicator based on the $P0$-projection~\eqref{eq:p0_proj_op}. 
The near-band~$\mathfrak{K}_h^\mathrm{near}$ is defined as the set of cells, which consists of the neighboring cells of all cut-cells.
Obviously, in the one-dimensional case, $\mathfrak{K}_h^\mathrm{near}$ consists of only two cells.
In Figure~\ref{fig:p0_indicator}, the exact density values~$\rho_\mathrm{pre}$ and~$\rho_\mathrm{post}$ of the 
pre- and post-shock region can be extracted from~$\mathfrak{K}_h^\mathrm{near}$. 
The width of the near-band may be extended in other test cases to extract~$\rho_\mathrm{pre}$ and~$\rho_\mathrm{post}$.

Subsequently, we define the local indicator~$\pzeroindicator_j (\rho)$ acting on the density~$\rho$ in a cut background cell~$K_j \in \mathfrak{K}_h$ with~$\cutcell \neq \emptyset$ as
\begin{align}	\label{eq:p0_indicator}
	\pzeroindicator_j (\rho)
	=
	\pzeroindicator (\rho) \vert_{K_j}
	\coloneqq
	\underbrace{\frac{\rho_\mathrm{pre} + \rho_\mathrm{post}}{2}}_{\substack{\textrm{exact solution in} \\ \textrm{near-band}~\mathfrak{K}_h^\mathrm{near}}}
	-
	\underbrace{\frac{\rho_{j, \mathfrak{A}}^\mathrm{P0} + \rho_{j, \mathfrak{B}}^\mathrm{P0}}{2}}_{\substack{P0\textrm{-projections} \\ \textrm{in cut-cells}~\cutcell}}
	\,,	\qquad
	\mathfrak{s} \in \{ \mathfrak{A}, \mathfrak{B} \}
	\,.
\end{align}
In Equation~\ref{eq:p0_indicator}, we compare the average value of the exact solution, which we extract from the near-band~$\mathfrak{K}_h^\mathrm{near}$, to the average value of the $P0$-projections in the cut-cells~$K_{j, \mathfrak{A}}$ and~$K_{j, \mathfrak{B}}$. Based on the sign of the indicator~$\pzeroindicator_j (\rho)$, the shifting direction of the shock interface~$\mathfrak{I}$ is determined by
\begin{subequations}	\label{eq:shifting_direction}
	\begin{align}	
		&\textrm{if } \pzeroindicator_j (\rho) > 0 \qquad \Rightarrow \qquad \textrm{shift }\mathfrak{I} \textrm{ to the right}\,,\\
		&\textrm{if } \pzeroindicator_j (\rho) < 0 \qquad \Rightarrow \qquad \textrm{shift }\mathfrak{I} \textrm{ to the left}\,.
	\end{align}
\end{subequations}
In other words, Equation~\ref{eq:shifting_direction} tells us the following: A value of~$\pzeroindicator_j (\rho) < 0 $ indicates that the average value of the $P0$-projections is larger than the exact average value of the shock. Thus, the shock interface~$\shockInterface$ is located too far right in relation to the exact shock position, i.e., $\interfacepos > \shockpos$.
We obtain the new shock interface position~$\interfacepos^{l+1}$ by a bisection algorithm, which takes the last two shock interface positions~$\interfacepos^{l}$ and~$\interfacepos^{l-1}$ into account. During the first iterations, the boundaries~$\partial K_j$ of the cut background cell are considered. We evaluate Equation~\eqref{eq:p0_indicator} for the scenario shown in Figure~\ref{fig:p0_indicator}, which gives an initial value of~$\pzeroindicator_j (\rho) = -0.389$. Bisection yields the new shock interface position~$\interfacepos^{l+1} = [x(\partial K_j^\mathrm{left}) + \interfacepos^l]/2 = [0.5 + 0.57] / 2 = 0.535$.

\subsubsection{Indicators Based on the Normal Shock Relations}	\label{sec:normal_shock_rel_indicators}
We state the one-dimensional normal shock relations~\cite{anderson1990}
\begin{subequations}	\label{eq:normalShockRelations}
	\begin{align}
		\rho_\mathrm{pre} u_{1,\mathrm{pre}} &= \rho_\mathrm{post} u_{1,\mathrm{post}} && \textrm{(continuity eq.)}\label{eq:normalShockWaveRelations_rpt_conti}\,,\\
		p_\mathrm{pre} +  \rho_\mathrm{pre} u_{1,\mathrm{pre}}^2 &= p_\mathrm{post} + \rho_\mathrm{post} u_{1,\mathrm{post}}^2 && \textrm{(momentum eq.)}\label{eq:normalShockWaveRelations_rpt_mom}\,,
	\end{align}
\end{subequations}
for the definition of two additional indicators. Equation~\eqref{eq:normalShockRelations} has to be fulfilled up to the accuracy of the numerical method if the shock is resolved sharply. Consequently, a deviation from this state motivates the definition of the indicators~$\rhoindicator (\vec{U})$ and~$\momindicator (\vec{U})$, which are based on Equations~\eqref{eq:normalShockWaveRelations_rpt_conti} and~\eqref{eq:normalShockWaveRelations_rpt_mom}, respectively.

We define the local indicator~$\rhoindicator_j (\vec{U})$ in a cut background cell~$K_j \in \mathfrak{K}_h$  with $\cutcell \neq \emptyset$ as
\begin{align}	\label{eq:rho_indicator}
	\begin{split}
		\rhoindicator_j (\vec{U})
		=
		\rhoindicator (\vec{U}) \vert_{K_j}
		&\coloneqq
		\frac{1}{\partial K_j}
		\oint_\mathfrak{I}
		\rho_{j, \mathfrak{A}} u_{j, \mathfrak{A}} - \rho_{j,\mathfrak{B}} 	u_{j, \mathfrak{B}}
		\diff{S}\\
		&=
		\frac{1}{\partial K_j}
		\oint_\mathfrak{I}
		\jump{\rho_j u_{j}}
		\diff{S}
		\,,	
	\end{split}
\end{align}
and the indicator~$\momindicator_j (\vec{U})$ as
\begin{equation}	\label{eq:indicator_mom}
	\begin{split}
		\momindicator_j (\vec{U})
		=
		\momindicator (\vec{U}) \vert_{K_j}
		&\coloneqq
		\frac{1}{\partial K_j}
		\oint_\mathfrak{I}
		p_{j, \mathfrak{A}} + \rho_{j,\mathfrak{A}} u_{j,\mathfrak{A}}^2 - \left(p_{j,\mathfrak{B}} + \rho_{j,\mathfrak{B}} u_{j,\mathfrak{B}}^2\right)
		\diff{S}\\
		&=
		\frac{1}{\partial K_j}
		\oint_\mathfrak{I}
		\left[\!\!\left[	
		p_j + \rho_j u_{j}^2
		\right]\!\!\right]
		\diff{S}
		\,.
	\end{split}
\end{equation}
The new shock interface position~$\interfacepos^{l+1}$ is computed by applying the bisection algorithm as presented for the indicator~$\pzeroindicator_j (\rho)$ in Section~\ref{sec:p0_indicator}.

\subsection{Implicit Pseudo-Time-Stepping}	\label{sec:implicit_pseudo_ts}
In the following, we present a novel implicit pseudo-time-stepping procedure in order to correct the shock interface position~$\interfacepos$ inside a cut background cell~$K_j$ by means of the indicators presented in Section~\ref{sec:indicators}. For that, we assume that a cell-accurate guess of~$\interfacepos$ is already known. In Section~\ref{sec:implicit_pseudo_ts_setting}, we explain the setting of the test case. We show the results of a one-dimensional proof of concept in Section~\ref{sec:implicit_pseudo_ts_proof}.

\subsubsection{Setting}	\label{sec:implicit_pseudo_ts_setting}
To illustrate the procedure, we consider a stationary shock wave located at~$\shockpos = 0.55$ with a Mach number of $\mathrm{M_s} = 1.5$ by fixing the frame of reference to it.
We choose a computational domain $\Omega_h = (0, 1)$ discretized by $10$ equally-sized cells with~$h=0.1$.
We prescribe 
the exact pre- and post-shock states at the left and right boundary, respectively. The pre-shock conditions are given by
\begin{align}	\label{eq:pre_shock_cond}
	\left( \rho_\mathrm{pre},\,u_{\mathrm{pre}},\,p_\mathrm{pre}\right)^\intercal
	=
	\left(1,\,\sqrt{\gamma \frac{p_\mathrm{pre}}{\rho_\mathrm{pre}}}\mathrm{M_s},\,1 \right)^\intercal
\end{align}
and the post-shock conditions by
\begin{subequations} \label{eq:post_shock_cond}
	\begin{align}
		\rho_\mathrm{post} &= \frac{(\gamma + 1) \mathrm{M_s}^2}{2 + (\gamma -1) \mathrm{M_s}^2} \rho_\mathrm{pre}\,,\\
	  	u_{\mathrm{post}} &= \frac{2 + (\gamma - 1) \mathrm{M_s}^2}{(\gamma + 1) \mathrm{M_s}^2} u_{1, \mathrm{pre}}\,,\\
		p_\mathrm{post} &= \left[ 1 + \frac{2 \gamma}{\gamma + 1} (\mathrm{M_s}^2 - 1) \right] p_\mathrm{pre}\,.
	\end{align}
\end{subequations}
In Equations~\eqref{eq:pre_shock_cond} and~$\eqref{eq:post_shock_cond}$, we set the heat capacity ratio to~$\gamma = 1.4$.

As a starting guess, we mimic a sufficiently smooth initial condition, which can be obtained by means of a smoothed Heaviside function
\begin{align}	\label{eq:smoothed_Heaviside}
	H (x)
	=
	0.5
	\left[
		\tanh \left( \frac{\normVector{x - x_\mathrm{s}}}{\tilde{C} h / \max(0, P)} + 1.0 \right)
	\right]
	\,,
\end{align}
where~$x$ denotes an arbitrary point in the domain, $x_\mathrm{s}$ is position of the shock front, 
and $\tilde{C} = 1.0$ is a user-defined factor to control the strength of the smoothing. 
The smoothed initial conditions of a physical quantity~$\psi_0 (x) = \psi (x, t_0)$ are then computed with
\begin{align}	\label{eq:smoothed_initial_cond}
	\psi_0 (x)=
	\psi_\mathrm{pre} (x) - H(x) \left[ \psi_\mathrm{pre} (x) - \psi_\mathrm{post} (x) \right]\,.
\end{align}
We apply the smoothing~\eqref{eq:smoothed_initial_cond} to all components of the state vector~$\vec{U}$, see Equation~\eqref{eq:vectorConservedQuantities}. Note that we do not expect any limitation in using a DG computation with a suitable shock-capturing strategy~\cite{geisenhofer2019,persson2006,barter2010} as an input for the XDG computation.

\begin{figure}[tbp]
	\centering
	\begin{tikzpicture}
		\begin{axis}[
			xmin=0.45, xmax=0.65,
			ymin=0.8, ymax=2.2,
			xlabel=$x$,
			ylabel=Density~$\rho$,
			xtick={0.5, 0.55, 0.6},
			xticklabels={$0.5$, $\shockpos$, $0.6$},
			no markers,
			xmajorgrids=false,
			ymajorgrids=false,
			legend style={cells={align=left}},
			width=0.27\textwidth,
			height=0.27\textwidth
			]
			
			\draw[color=lightgray, thick] (axis cs:0.5,\pgfkeysvalueof{/pgfplots/ymin}) -- (axis cs:0.5,\pgfkeysvalueof{/pgfplots/ymax});	
			\draw[color=lightgray, thick] (axis cs:0.6,\pgfkeysvalueof{/pgfplots/ymin}) -- (axis cs:0.6,\pgfkeysvalueof{/pgfplots/ymax});
			\draw[color=lightgray, dashed] (axis cs:0.55,\pgfkeysvalueof{/pgfplots/ymin}) -- (axis cs:0.55,\pgfkeysvalueof{/pgfplots/ymax});		
			
			\draw[color=red] (axis cs:0.57,0.8) -- (axis cs:0.57,2.2);
			\node[color=red, align=left, anchor=south west] at (axis cs:0.57,0.8) {\small$\mathfrak{I}_0$};			
			
			\pgfplotstableread{XDG_SSW_Initial_Condition.txt}\datafile
			\addplot+[] table[skip first n=1] {\datafile};
			\addlegendentry{Shock};	
			
			\pgfplotstableread{XDG_SSW_OneLs_p2_xCells10_yCells3_agg0.3_ts=100_dtFixed0.1_shockPos=0.55_smooth=1.txt}\datafile
			\addplot+[] table[skip first n=1] {\datafile};
			\addlegendentry{Smoothed initial\\condition, $P=2$};
			
		\end{axis}	
	\end{tikzpicture}
	\caption
	{Initial configuration for determining the sub-cell accurate position of a stationary shock wave in an XDG method. The smoothed initial condition mimics the solution of standard DG method with shock-capturing, see Equations~\eqref{eq:smoothed_Heaviside} and~\eqref{eq:smoothed_initial_cond}. The shock is located at~$\shockpos = 0.55$. The initial position of the shock interface~$\mathfrak{I}_0$ is assumed to be at~$x_\mathfrak{I}^0 = \interfacepos (t_0) = 0.57$.}
	\label{fig:initial_cond}
\end{figure}
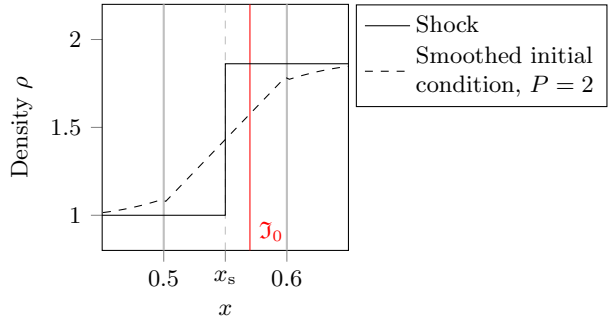

In order to verify the robustness of the sub-cell correction procedure, we assume the shock interface~$\mathfrak{I}_0$ (we skip the index~$\mathrm{s}$ which denotes the interface as the \emph{shock interface} in this section in order avoid confusion) to be initially located at~$x_\mathfrak{I}^{l=0} = x_\mathfrak{I}^0 = \interfacepos (t_0) = 0.57$. This position is sufficiently far away from the exact shock position~$\shockpos = 0.55$. The shock interface is represented by the zero iso-contour of a level-set function which is given by
\begin{align}	\label{eq:level_set}
	\shockLevelSet (x)
	=
	x - x_\mathfrak{I}^0
	=
	x - 0.57\,.
\end{align}
We show the smoothed initial condition centered around the shock interface~$\mathfrak{I}_0$ in Figure~\ref{fig:initial_cond}.
Note that the shock interface is fixed in space in each correction step. The procedure terminates if the position~$\interfacepos^l$ of the shock interface~$\mathfrak{I}_l$ coincides with the exact shock position~$\shockpos$, i.e., $\interfacepos^l = \shockpos$.

In order to obtain the steady-state solution during each correction step~$t_l$ of the shock interface position,
which we call a \emph{pseudo-time-step}, 
we use a standard implicit Euler scheme with finite time-step sizes of ~$\Delta t_l = 0.1$. 
This mathematically corresponds to an under-relaxation-like method, 
and thus the nonlinear system can be solved using a standard Newton procedure.
Due to the low number of DOFs the Jacobian matrix in some point can be evaluated in a brute-force fashion using a 
small perturbation in the order of $10^{-7}$ of each DOF.
The software framework \emph{BoSSS} provides some acceleration of this process, which exploits the locality of the DG discretization:
Since any perturbation only affects the residual in the respective cell and all its neighboring cells, one can perturb multiple DOFs at once,
given that they are sufficiently far apart. This allows to construct several columns of the Jacobian matrix at once. 
Further advanced preconditioning techniques and implicit methods~\cite{dolejsi2004,fidkowski2005,shahbazi2009} are beyond the scope of this work, 
since we focus on the novel methodology for determining sub-cell accurate position of the shock front.

We start the correction procedure from an initial numerical solution~$\vec{U}_0 (x) = \vec{U} (x, t_0)$ 
and an initial cell-accurate approximation of the shock interface position~$x_\mathfrak{I}^0$. 
Next, we set up the iterative pseudo-time-stepping procedure~$t_l \rightarrow t_{l+1}$ 
by discretizing the Euler equations~\eqref{eq:eulerEquations} to~\eqref{eq:convectiveFluxes} 
in combination with the shock interface~$\mathfrak{I}_l$
by the presented XDG method in space, see Section~\ref{sec:discretization}.
For obtaining the steady-state solution~$\vec{U} (x, t_l)$ of each pseudo-time-step~$t_l$, 
we apply the implicit Euler scheme to advance the solution in time~$t_l^n \rightarrow t^{n+1}_l$. 
Note that the index~$l$ corresponds to a pseudo-time-step where the shock interface~$\mathfrak{I}_l$ is fixed in space and the index~$n$ denotes the time-step of the implicit Euler scheme.

\subsubsection{Proof of Concept} \label{sec:implicit_pseudo_ts_proof}
\begin{figure}[tbp]
	\centering
	\begin{tikzpicture}
		\begin{groupplot}[
			group style={
				group name=mygroupplot,
				group size=3 by 4,
				vertical sep = 10mm,
				horizontal sep = 19mm,
			},
			height=0.2\textwidth,
			width=0.2\textwidth,
			no markers
			]
			
			% Plot 1
			\nextgroupplot[
			xmin=0.45, xmax=0.65,
			ymin=0, ymax=1,
            ylabel={Initial value},
			xmajorgrids=false,
			ymajorgrids=false,
			yticklabels={,,},
			ytick style={draw=none},
			xtick={0.5, 0.55, 0.6},
			xticklabels={$0.5$, $\shockpos$, $0.6$},			
			title style={at={(0.5,1.0)}, anchor=north, yshift=5mm},
			title ={\small Cut background cell~$K_j$},
			]
			\draw[color=lightgray, thick] (axis cs:0.5,\pgfkeysvalueof{/pgfplots/ymin}) -- (axis cs:0.5,\pgfkeysvalueof{/pgfplots/ymax});	
			\draw[color=lightgray, thick] (axis cs:0.6,\pgfkeysvalueof{/pgfplots/ymin}) -- (axis cs:0.6,\pgfkeysvalueof{/pgfplots/ymax});
			\draw[color=lightgray, dashed] (axis cs:0.55,\pgfkeysvalueof{/pgfplots/ymin}) -- (axis cs:0.55,\pgfkeysvalueof{/pgfplots/ymax});					
			\draw[color=red, thick] (axis cs:0.57,\pgfkeysvalueof{/pgfplots/ymin}) -- (axis cs:0.57,\pgfkeysvalueof{/pgfplots/ymax});			
			\node[color=red, align=left, anchor=south west] at (axis cs:0.57,0) {\small $\mathfrak{I}_0$ at\\\small $0.57$};			
			\node[color=black,fill=black!10!,anchor=north east] at (axis cs:\pgfkeysvalueof{/pgfplots/xmax},\pgfkeysvalueof{/pgfplots/ymax}) {\small$t_0$};
			
			% Plot 2
			\nextgroupplot[
			xmin=0.45, xmax=0.65,
			ymin=0.75, ymax=2.15,
			ylabel=Density~$\rho$,
			xmajorgrids=false,
			ymajorgrids=false,
			no markers,
			xtick={0.5, 0.55, 0.6},
			xticklabels={$0.5$, $\shockpos$, $0.6$},
			legend columns=-1,
			column sep=0.1cm,
			legend to name={CommonLegend},
			title style={at={(0.5,1.0)}, anchor=north, yshift=5mm},
			title ={\small Numerical solution},
			]
			\draw[color=lightgray, thick] (axis cs:0.5,\pgfkeysvalueof{/pgfplots/ymin}) -- (axis cs:0.5,\pgfkeysvalueof{/pgfplots/ymax});	
			\draw[color=lightgray, thick] (axis cs:0.6,\pgfkeysvalueof{/pgfplots/ymin}) -- (axis cs:0.6,\pgfkeysvalueof{/pgfplots/ymax});								
			\node[color=black,fill=black!10!,anchor=north east] at (axis cs:\pgfkeysvalueof{/pgfplots/xmax},\pgfkeysvalueof{/pgfplots/ymax}) {\small$t_0$};	
			
			\pgfplotstableread{XDG_SSW_Initial_Condition.txt}\datafile
			\addplot+[color=lightgray, dashed] table[skip first n=1] {\datafile};
			
			\pgfplotstableread{1_rho.txt}\datafile
			\addplot+[solid] table[skip first n=1] {\datafile};
			
			\pgfplotstableread{1_rho_p0.txt}\datafile
			\addplot+[dashdotted] table[skip first n=1] {\datafile};
			
			\legend{Shock,$P=2$, $P0$-proj.}
			
			% Plot 3
			\nextgroupplot[
			xmin=-0.5, xmax=2.5,
			ymin=-1.0, ymax=5.0,
			ylabel=Indicator value,
			xmajorgrids=true,
			ymajorgrids=true,
			xtick={0, 1, 2},
			xticklabels={$\pzeroindicator$, $\rhoindicator$, $\momindicator$},
			title style={at={(0.5,1.0)}, anchor=north, yshift=5mm},
			title ={\small Indicators},
			]
			\pgfplotstableread{1_indicators_p0.txt}\datafile
			\addplot+[
			ybar,
			fill=gray,
			] table[] {\datafile};
			
			\pgfplotstableread{1_indicators.txt}\datafile
			\addplot+[
			ybar,
			fill=black!10!,
			solid
			] table[] {\datafile};
			
			\node[color=black,fill=black!10!,anchor=north east] at (axis cs:\pgfkeysvalueof{/pgfplots/xmax},\pgfkeysvalueof{/pgfplots/ymax}) {\small $t_0$};	
			
			% Plot 4
			\nextgroupplot[
			xmin=0.45, xmax=0.65,
			ymin=0, ymax=1,
            ylabel={Time-step 1},
			xmajorgrids=false,
			ymajorgrids=false,
			yticklabels={,,},
			ytick style={draw=none},
			xtick={0.5, 0.55, 0.6},
			xticklabels={$0.5$, $\shockpos$, $0.6$},			
			]
			\draw[color=lightgray, thick] (axis cs:0.5,\pgfkeysvalueof{/pgfplots/ymin}) -- (axis cs:0.5,\pgfkeysvalueof{/pgfplots/ymax});	
			\draw[color=lightgray, thick] (axis cs:0.6,\pgfkeysvalueof{/pgfplots/ymin}) -- (axis cs:0.6,\pgfkeysvalueof{/pgfplots/ymax});
			\draw[color=lightgray, dashed] (axis cs:0.55,\pgfkeysvalueof{/pgfplots/ymin}) -- (axis cs:0.55,\pgfkeysvalueof{/pgfplots/ymax});					
			\draw[color=red, thick] (axis cs:0.535,\pgfkeysvalueof{/pgfplots/ymin}) -- (axis cs:0.535,\pgfkeysvalueof{/pgfplots/ymax});			
			\node[color=red, align=left, anchor=south west] at (axis cs:0.535,0) {\small$\mathfrak{I}_1$ at\\\small$0.535$};			
			\node[color=black,fill=black!10!,anchor=north east] at (axis cs:\pgfkeysvalueof{/pgfplots/xmax},\pgfkeysvalueof{/pgfplots/ymax}) {\small$t_1$};
			
			% Plot 5			
			\nextgroupplot[
			xmin=0.45, xmax=0.65,
			ymin=0.75, ymax=2.15,
			ylabel=Density~$\rho$,
			xmajorgrids=false,
			ymajorgrids=false,
			no markers,
			xtick={0.5, 0.55, 0.6},
			xticklabels={$0.5$, $\shockpos$, $0.6$},
			legend columns=-1,
			column sep=0.1cm,
			legend to name={CommonLegend}
			]
			\draw[color=lightgray, thick] (axis cs:0.5,\pgfkeysvalueof{/pgfplots/ymin}) -- (axis cs:0.5,\pgfkeysvalueof{/pgfplots/ymax});	
			\draw[color=lightgray, thick] (axis cs:0.6,\pgfkeysvalueof{/pgfplots/ymin}) -- (axis cs:0.6,\pgfkeysvalueof{/pgfplots/ymax});								
			\node[color=black,fill=black!10!,anchor=north east] at (axis cs:\pgfkeysvalueof{/pgfplots/xmax},\pgfkeysvalueof{/pgfplots/ymax}) {\small$t_1$};	
			
			\pgfplotstableread{XDG_SSW_Initial_Condition.txt}\datafile
			\addplot+[color=lightgray, dashed] table[skip first n=1] {\datafile};
			
			\pgfplotstableread{2_rho.txt}\datafile
			\addplot+[solid] table[skip first n=1] {\datafile};
			
			\pgfplotstableread{2_rho_p0.txt}\datafile
			\addplot+[dashdotted] table[skip first n=1] {\datafile};
			
			% Plot 6			
			\nextgroupplot[
			xmin=-0.5, xmax=2.5,
			ymin=-1.0, ymax=5.0,
			ylabel=Indicator value,
			xmajorgrids=true,
			ymajorgrids=true,
			xtick={0, 1, 2},
			xticklabels={$\pzeroindicator$, $\rhoindicator$, $\momindicator$},
			]
			\pgfplotstableread{2_indicators_p0.txt}\datafile
			\addplot+[
			ybar,
			fill=gray,
			] table[] {\datafile};
			
			\pgfplotstableread{2_indicators.txt}\datafile
			\addplot+[
			ybar,
			fill=black!10!,
			solid
			] table[] {\datafile};
			
			\node[color=black,fill=black!10!,anchor=north east] at (axis cs:\pgfkeysvalueof{/pgfplots/xmax},\pgfkeysvalueof{/pgfplots/ymax}) {\small$t_1$};	
			
			% Plot 7			
			\nextgroupplot[
			xmin=0.45, xmax=0.65,
			ymin=0, ymax=1,
            ylabel={Time-step 2},
			xmajorgrids=false,
			ymajorgrids=false,
			yticklabels={,,},
			ytick style={draw=none},
			xtick={0.5, 0.55, 0.6},
			xticklabels={$0.5$, $\shockpos$, $0.6$},			
			]
			\draw[color=lightgray, thick] (axis cs:0.5,\pgfkeysvalueof{/pgfplots/ymin}) -- (axis cs:0.5,\pgfkeysvalueof{/pgfplots/ymax});	
			\draw[color=lightgray, thick] (axis cs:0.6,\pgfkeysvalueof{/pgfplots/ymin}) -- (axis cs:0.6,\pgfkeysvalueof{/pgfplots/ymax});
			\draw[color=lightgray, dashed] (axis cs:0.55,\pgfkeysvalueof{/pgfplots/ymin}) -- (axis cs:0.55,\pgfkeysvalueof{/pgfplots/ymax});					
			\draw[color=red, thick] (axis cs:0.5525,\pgfkeysvalueof{/pgfplots/ymin}) -- (axis cs:0.5525,\pgfkeysvalueof{/pgfplots/ymax});			
			\node[color=red, align=left, anchor=south west] at (axis cs:0.5525,0) {\small$\mathfrak{I}_2$ at\\\small$0.5525$};			
			\node[color=black,fill=black!10!,anchor=north east] at (axis cs:\pgfkeysvalueof{/pgfplots/xmax},\pgfkeysvalueof{/pgfplots/ymax}) {\small$t_2$};
			
			% Plot 8			
			\nextgroupplot[
			xmin=0.45, xmax=0.65,
			ymin=0.75, ymax=2.15,
			ylabel=Density~$\rho$,
			xmajorgrids=false,
			ymajorgrids=false,
			no markers,
			xtick={0.5, 0.55, 0.6},
			xticklabels={$0.5$, $\shockpos$, $0.6$},
			legend columns=-1,
			column sep=0.1cm,
			legend to name={CommonLegend}
			]
			\draw[color=lightgray, thick] (axis cs:0.5,\pgfkeysvalueof{/pgfplots/ymin}) -- (axis cs:0.5,\pgfkeysvalueof{/pgfplots/ymax});	
			\draw[color=lightgray, thick] (axis cs:0.6,\pgfkeysvalueof{/pgfplots/ymin}) -- (axis cs:0.6,\pgfkeysvalueof{/pgfplots/ymax});								
			\node[color=black,fill=black!10!,anchor=north east] at (axis cs:\pgfkeysvalueof{/pgfplots/xmax},\pgfkeysvalueof{/pgfplots/ymax}) {\small$t_2$};	
			
			\pgfplotstableread{XDG_SSW_Initial_Condition.txt}\datafile
			\addplot+[color=lightgray, dashed] table[skip first n=1] {\datafile};
			
			\pgfplotstableread{3_rho.txt}\datafile
			\addplot+[solid] table[skip first n=1] {\datafile};
			
			\pgfplotstableread{3_rho_p0.txt}\datafile
			\addplot+[dashdotted] table[skip first n=1] {\datafile};
			
			% Plot 9			
			\nextgroupplot[
			xmin=-0.5, xmax=2.5,
			ymin=-1.0, ymax=5.0,
			ylabel=Indicator value,
			xmajorgrids=true,
			ymajorgrids=true,
			xtick={0, 1, 2},
			xticklabels={$\pzeroindicator$, $\rhoindicator$, $\momindicator$},
			]
			\pgfplotstableread{3_indicators_p0.txt}\datafile
			\addplot+[
			ybar,
			fill=gray,
			] table[] {\datafile};
			
			\pgfplotstableread{3_indicators.txt}\datafile
			\addplot+[
			ybar,
			fill=black!10!,
			solid
			] table[] {\datafile};
			
			\node[color=black,fill=black!10!,anchor=north east] at (axis cs:\pgfkeysvalueof{/pgfplots/xmax},\pgfkeysvalueof{/pgfplots/ymax}) {\small$t_2$};	
			
			% Plot 10			
			\nextgroupplot[
			xmin=0.45, xmax=0.65,
			ymin=0, ymax=1,
			xlabel=$x$,
            ylabel={Final time-step},
			xmajorgrids=false,
			ymajorgrids=false,
			yticklabels={,,},
			ytick style={draw=none},
			xtick={0.5, 0.55, 0.6},
			xticklabels={$0.5$, $\shockpos$, $0.6$},			
			]
			\draw[color=lightgray, thick] (axis cs:0.5,\pgfkeysvalueof{/pgfplots/ymin}) -- (axis cs:0.5,\pgfkeysvalueof{/pgfplots/ymax});	
			\draw[color=lightgray, thick] (axis cs:0.6,\pgfkeysvalueof{/pgfplots/ymin}) -- (axis cs:0.6,\pgfkeysvalueof{/pgfplots/ymax});
			\draw[color=lightgray, dashed] (axis cs:0.55,\pgfkeysvalueof{/pgfplots/ymin}) -- (axis cs:0.55,\pgfkeysvalueof{/pgfplots/ymax});					
			\draw[color=red, thick] (axis cs:0.55,\pgfkeysvalueof{/pgfplots/ymin}) -- (axis cs:0.55,\pgfkeysvalueof{/pgfplots/ymax});			
			\node[color=red, align=left, anchor=south west] at (axis cs:0.55,0) {\small$\mathfrak{I}_\mathrm{last}$ at\\\small$0.55$};			
			\node[color=black,fill=black!10!,anchor=north east] at (axis cs:\pgfkeysvalueof{/pgfplots/xmax},\pgfkeysvalueof{/pgfplots/ymax}) {\small$t_\mathrm{last}$};
			
			% Plot 11			
			\nextgroupplot[
			xmin=0.45, xmax=0.65,
			ymin=0.75, ymax=2.15,
			xlabel=$x$,
			ylabel=Density~$\rho$,
			xmajorgrids=false,
			ymajorgrids=false,
			no markers,
			xtick={0.5, 0.55, 0.6},
			xticklabels={$0.5$, $\shockpos$, $0.6$},
			legend columns=-1,
			legend style={/tikz/every even column/.append style={column sep=3mm}},
			legend to name={CommonLegend}
			]
			\draw[color=lightgray, thick] (axis cs:0.5,\pgfkeysvalueof{/pgfplots/ymin}) -- (axis cs:0.5,\pgfkeysvalueof{/pgfplots/ymax});	
			\draw[color=lightgray, thick] (axis cs:0.6,\pgfkeysvalueof{/pgfplots/ymin}) -- (axis cs:0.6,\pgfkeysvalueof{/pgfplots/ymax});							
			\node[color=black,fill=black!10!,anchor=north east] at (axis cs:\pgfkeysvalueof{/pgfplots/xmax},\pgfkeysvalueof{/pgfplots/ymax}) {\small$t_\mathrm{last}$};	
			
			\pgfplotstableread{XDG_SSW_Initial_Condition.txt}\datafile
			\addplot+[color=lightgray, dashed] table[skip first n=1] {\datafile};
			
			\pgfplotstableread{5_rho.txt}\datafile
			\addplot+[solid] table[skip first n=1] {\datafile};
			
			\pgfplotstableread{5_rho_p0.txt}\datafile
			\addplot+[dashdotted] table[skip first n=1] {\datafile};
			
			\legend{Shock, {$P=2$}, {$P0$-proj.}}
			
			% Plot 12			
			\nextgroupplot[
			xmin=-0.5, xmax=2.5,
			ymin=-1.0, ymax=5.0,
			xlabel style={align=center},
			xlabel={Indicator variant\\$\pzeroindicator$: $+/- \Rightarrow$\\$\mathfrak{I} \textrm{: to right/left}$},
			ylabel=Indicator value,
			xmajorgrids=true,
			ymajorgrids=true,
			xtick={0, 1, 2},
			xticklabels={$\pzeroindicator$, $\rhoindicator$, $\momindicator$},
			]
			\pgfplotstableread{5_indicators_p0.txt}\datafile
			\addplot+[
			ybar,
			fill=gray,
			] table[] {\datafile};
			
			\pgfplotstableread{5_indicators.txt}\datafile
			\addplot+[
			ybar,
			fill=black!10!,
			solid
			] table[] {\datafile};
			
			\node[color=black,fill=black!10!,anchor=north east] at (axis cs:\pgfkeysvalueof{/pgfplots/xmax},\pgfkeysvalueof{/pgfplots/ymax}) {\small$t_\mathrm{last}$};				
		\end{groupplot}			
		\path ([yshift=-10mm]mygroupplot c1r4.south east) -- node[below]{\ref{CommonLegend}} ([yshift=-10mm]mygroupplot c3r4.south west);
	\end{tikzpicture}
\caption{
Illustration of the implicit pseudo-time-stepping procedure. 
The shock interface~$\mathfrak{I}_l$ converges to the exact shock position~$\shockpos=0.55$ by using the indicator~$\pzeroindicator (\rho)$. 
Three different indicators are presented: $\pzeroindicator (\rho)$ is based on a $P0$-projection of a high-order solution ($P=2$) of the density field~$\rho$, 
whereas $\rhoindicator (\vec{U})$ and $\momindicator (\vec{U})$ are based on 
the normal shock relations~\eqref{eq:normalShockWaveRelations_rpt_conti} and~\eqref{eq:normalShockWaveRelations_rpt_mom}, respectively. 
The shock interface position inside the cut background cell~$K_j = (0.5, 0.6)$ is iteratively adapted by applying a bisection procedure 
depending on the sign of~$\pzeroindicator (\rho)$. 
In this test case, only the indicator~$\pzeroindicator (\rho)$ delivers the correct shifting direction. 
The pseudo-time-steps~$t_l = \{t_0,t_1,t_2,t_\mathrm{last}\}$ are shown with gray labels.
}
\label{fig:proof_of_concept}
\end{figure}
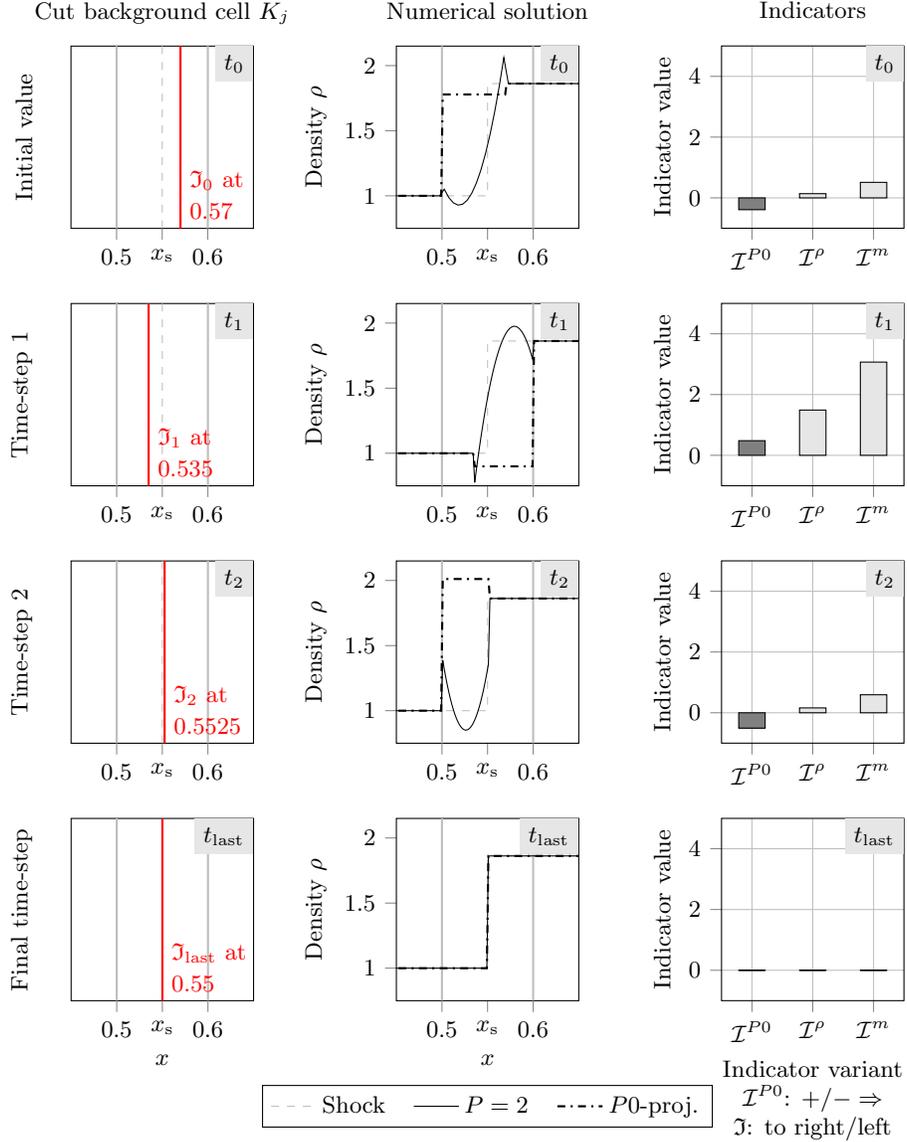

In Figure~\ref{fig:proof_of_concept}, we show the results of a one-dimensional proof of concept of the sub-cell accurate correction procedure in a pseudo-two-dimensional computation.
The reader is referred to Section~\ref{sec:implicit_pseudo_ts_setting} for a description of the test case.
The left column depicts the position of the shock interface~$\mathfrak{I}_l$ in the cut background cell~$K_j = (0.5, 0.6)$,
the middle column shows the unlimited solution with a polynomial degree of~$P=2$ as well as the corresponding $P0$-projections in the cut-cells~$K_{j, \mathfrak{A}}$ and~$K_{j, \mathfrak{B}}$, 
and the right column depicts the values of the indicators~$\pzeroindicator (\rho)$, $\rhoindicator (\vec{U})$, and $\momindicator (\vec{U})$. 
The results are plotted at the pseudo-time-steps~$t_l = \{t_0,t_1,t_2\}$, 
and for the converged solution at~$t_\mathrm{last}$ where the shock interface coincides with the exact shock position, i.e., $x_{\mathfrak{I}}^\mathrm{last} = \shockpos = 0.55$.

During the correction, only the indicator~$\pzeroindicator (\rho)$ is capable to deliver the correct shifting direction of the shock interface in every pseudo-time-step. All indicators give the correct shifting direction at the pseudo-time-step~$t_1$. By contrast, at the pseudo-time-steps~$t_0$ and~$t_2$, the indicators~$\rhoindicator (\vec{U})$ and~$\momindicator (\vec{U})$ based on the normal shock relations~\eqref{eq:normalShockRelations} fail to return the correct shifting direction, see again the right column of Figure~\ref{fig:proof_of_concept}. The new interface positions~$\interfacepos^{l+1}$ are determined by the bisection algorithm as discussed for the indicator~$\pzeroindicator (\rho)$ in Section~\ref{sec:p0_indicator}.

% ###############################################################################
% ###############################################################################
\section{Conclusion}	
\label{sec:conclusion}
% ###############################################################################
% ###############################################################################
In this work, we have derived a novel reconstruction procedure for determining the sub-cell accurate position of a pseudo-two-dimensional shock front in the context of shock-fitting by employing an XDG method. The overall goal is the simulation of compressible flows with high-order accuracy, for which the presented methodology builds a fundamental basis. 

Many DG shock-capturing strategies lack an (arbitrary) high-order convergence rate. By contrast, we employ an XDG method in which a shock front is described by the zero iso-contour of a level-set function in order to circumvent this limitation. To this end, it is essential to obtain a sharp reconstruction of the shock front. This is a non-trivial task, since the position and shape of shock waves are in general not known a priori. We have introduced the methodology to determine the exact position of a stationary normal shock wave in a pseudo-two-dimensional setting. Furthermore, our approach features a cell-agglomeration strategy~\cite{muller2017,kummer2020} for the treatment of small and ill-shaped cut-cells.

In Section~\ref{sec:sub-cell_accurate_correction}, we have presented a novel \emph{sub-cell accurate} correction procedure in order to find the \emph{exact} shock position in a cut background cell.
We have developed an implicit pseudo-time-stepping procedure which determines the exact shock-position in an iterative way. In each pseudo-time-step, we use a standard implicit Euler scheme to drive the solution into the steady-state.
In every pseudo-time-step, the position of the shock front in the cut background cell is judged by several local indicators. In this test case, the indicator based on a zeroth-order projection of the local solution has always determined the correct shifting direction of the shock interface. The new shock interface positions have been computed by using a bisection algorithm.

Our ongoing and future work focuses on the extension to truly two-di\-men\-sion\-al, steady applications. For that, several marker points may be seeded on the reconstructed shock interface. An adapted implicit pseudo time-stepping procedure could be applied to the marker points until they all converge to the exact position on the shock front. Then, the shock level-set function could be reconstructed.
For unsteady flow scenarios, the complexity may increase drastically, since at this point it is not foreseeable whether additional stabilization mechanism are required for the tracking of moving shock fronts.

% ###############################################################################
% ###############################################################################
\section*{Acknowledgments}
% ###############################################################################
% ###############################################################################
This work is supported by the \emph{Excellence Initiative} of the German Federal and State Governments and the Graduate School CE within the Center of Computational Engineering at the Technical University of Darmstadt.
Thanks to L. Beck and J. Guti\'{e}rrez for the valuable comments on this work.

% ###############################################################################
% ###############################################################################
\bibliographystyle{plain}
\bibliography{bibliography}
% ###############################################################################
% ###############################################################################

\end{document}